\documentclass[11pt,reqno]{article}
\usepackage{bbding}
\usepackage{pifont}
\usepackage{amsmath}
\usepackage{mathrsfs}
\usepackage{amsfonts,bm}
\usepackage{amsthm}
\usepackage{epsfig}
\usepackage[]{harpoon}
\usepackage[active]{srcltx}
\usepackage{indentfirst,latexsym}
\usepackage{psfrag}
\usepackage[]{amssymb}
\usepackage{cite}
\usepackage{caption}

\newcommand{\new}{\newcommand*}\new{\rnew}{\renewcommand*}
\new{\newe}{\newenvironment*}\new{\stl}{\setlength}
\stl{\textwidth}{155mm}\stl{\textheight}{22cm}\stl{\headheight}{0cm}
\stl{\topmargin}{0cm}\stl{\oddsidemargin}{0.5cm}\stl{\evensidemargin}{0cm}
\rnew{\arraystretch}{1.2}\rnew{\baselinestretch}{1.2}
\renewcommand{\thefootnote}{\ding{73}}
\newtheorem{thm}{Theorem}

\newcommand{\eps}{\varepsilon}

\newcommand{\dps}{\displaystyle}

\newcommand{\fr}{\frac}

\newcommand{\pa}{\partial}

\numberwithin{equation}{section}
\newe{keywords}
   {\begin{quote}{\bf Keywords:}}
      {\end{quote}}
\newe{AMS}
   {\begin{quote}{\bf AMS subject classification 2010:}}
      {\end{quote}}
\newe{MSC}{\vspace*{5mm}
    {\noindent\it Mathematics Subject Classification(2010):}}{}
\new{\sect}[1]{\section{#1}\setcounter{equation}{0}
 \setcounter{thm}{0}\setcounter{lmm}{0}\setcounter{rmk}{0} }

\begin{document}

\title{ Singularity for a multidimensional variational wave equation
arising from nematic liquid crystals}

\author{
Yanbo Hu$^{a,*}$, Guodong Wang$^{b}$
\\{\small \it $^a$Department of Mathematics, Hangzhou Normal University,
Hangzhou, 311121, PR China}
\\
{\small \it $^b$School of Mathematics \& Physics, Anhui Jianzhu University, Hefei, 230601,
PR China}}

\rnew{\thefootnote}{\fnsymbol{footnote}}

\footnotetext{ $^*$Corresponding author. }
\footnotetext{ Email address: yanbo.hu@hotmail.com (Y. Hu), yxgdwang@163.com (G. Wang). }

\date{}

\maketitle
\begin{abstract}
This article is focused on a multidimensional nonlinear variational wave equation which is the Euler-Lagrange equation of a variational principle arising form the theory of nematic liquid crystals. By using the method of characteristics, we show that the smooth solutions for the spherically-symmetric variational wave equation breakdown in finite time, even for the arbitrarily small initial energy.
\end{abstract}

\begin{keywords}
Variational wave equation; singularity; characteristic method
\end{keywords}

\begin{AMS}
35L05; 35L72; 35B44
\end{AMS}

\section{Introduction}\label{S1}

We are interested in a nonlinear wave equation derivable form a variational principle of
the theory of nematic liquid crystals. In nematic liquid crystals, the mean orientation of the long molecules
can be described by the so-called director field $\textbf{n}(\textbf{x},t)$ at a spatial location $\textbf{x}$ and time $t$. In the regime in which inertia effects dominate viscosity, that is, ignoring
the kinetic energy of the director field, the propagation of the orientation waves
in the director field can be modeled by the least action principle \cite{Hunter-Saxton}
\begin{align}\label{1.1}
\delta\int \bigg(\fr{1}{2}\pa_t\textbf{n}\cdot\pa_t\textbf{n}-W(\textbf{n},\nabla \textbf{n})\bigg){\rm d}\textbf{x}{\rm d}t=0,\qquad \textbf{n}\cdot\textbf{n}=1,
\end{align}
where $W(\textbf{n},\nabla \textbf{n})$ is the well-known Oseen-Franck potential energy density
\begin{align}\label{1.2}
W(\textbf{n},\nabla \textbf{n})=\fr{1}{2}k_1(\nabla\cdot\textbf{n})^2 +\fr{1}{2}k_2(\textbf{n}\cdot\nabla\times\textbf{n})^2 +\fr{1}{2}k_3|\textbf{n}\times(\nabla\times\textbf{n})|^2,
\end{align}
with elastic constants $k_1, k_2$ and $k_3$ of the material. For the planar deformations of the director field $\textbf{n}$ depending only on a single spatial variable with $\textbf{x}=(x, 0, 0)$ and $\textbf{n}=(\cos u(x,t), \sin u(x,t), 0)$, then \eqref{1.1} reduces to
\begin{align}\label{1.3}
\delta\int \bigg(\fr{1}{2}u_{t}^2-\fr{1}{2}c^2(u)u_{x}^2\bigg){\rm d}x{\rm d}t=0,
\end{align}
with $c^2(u)=k_1\sin^2u+k_3\cos^2u$. The Euler-Lagrange equation of \eqref{1.3} given by
\begin{align}\label{1.4}
u_{tt}-c(u)(c(u)u_x)_x=0,
\end{align}
which is also called the one-dimensional variational wave equation, see \cite{Ali, Hunter-Saxton} for more information and applications of \eqref{1.4}.

Equation \eqref{1.3} has been widely explored since it was introduced. In particular, Glassey, Hunter and Zheng \cite{Glassey1} shown that, even for smooth initial data with arbitrarily small energy, its solutions can form cusp-type singularities in finite time. The existence of global dissipative weak solutions to \eqref{1.4} were systematically studied by Zhang and Zheng in \cite{Zha-Zhe4, Zha-Zhe6, Zha-Zhe7}, also see the work of Bressan and Huang \cite{Bres-Huang}. In \cite{Bres-Zheng}, Bressan and Zheng proposed the the method of energy-dependent coordinates to establish the global existence of conservative solutions to its Cauchy problem for initial data of finite energy. The uniqueness of the conservative solutions was provided by Bressan, Chen and Zhang \cite{B-C-Z}. For more relevant results of \eqref{1.4} can be found among others in \cite{B-C, B-C1, Hu3, Hu-Wang}. We also refer the reader to Refs. \cite{Cai, Chen-Zheng, Zha-Zhe, Zha-Zhe9} for the study of the one-dimensional nonlinear variational wave systems.

The multidimensional version of variational wave equation \eqref{1.1} reads that \cite{Bres-Zheng, Glassey1, Zha-Zhe4}
\begin{align}\label{1.5}
u_{tt}-c(u)\nabla\cdot(c(u)\nabla u)=0.
\end{align}
Equation \eqref{1.5} can also be derived from the variational principle whose action is a quadratic function of the derivatives of the field \cite{Ali0, Ali, Hunter-Saxton}. In this short article, we show the singularity formation of smooth solutions for the nonlinear variational wave equation \eqref{1.5} with spherically symmetric initial data.
Set $u=u(t,r), r=|\textbf{x}|$, where $\textbf{x}=(x_1,\cdots,x_d)\ (d>1)$ are the spatial independent variables, then \eqref{1.5} can be reduced to
\begin{align}\label{1.6}
u_{tt}-c(u)(c(u) u_r)_r-\fr{(d-1)c^2(u)u_r}{r}=0.
\end{align}
We assume that the wave speed $c(\cdot)\in C^2$ satisfies
\begin{align}\label{1.7}
0<c_0\leq c(\cdot)\leq c_1,\quad |c'(\cdot)|\leq c_1
\end{align}
for some positive constants $c_0,c_1$. Our main singularity formation result is
\begin{thm}\label{thm}
Assume that $c'(u_0)>0$ for some constant $u_0$. Let $\phi(z)$ be a smooth function satisfying
\begin{align}\label{1.8}
\phi(z)\in C^{1}_c((-1,1)),\quad \phi'(0)\leq -2\max\bigg\{\fr{32c_{1}^22^{\alpha}}{r_0c_0c'(u_0)},\ \fr{1}{c_0r_{0}^{\alpha}}\bigg\},
\end{align}
where $r_0>0$ and $\alpha=(d-1)/2$. Suppose that $u(t,r)\in C^1$ is a smooth solution of \eqref{1.6} in $0\leq t<T, r>0$ with initial data
\begin{align}\label{1.9}
u(0,r)=u_0+\eps\phi\bigg(\dps\fr{r-r_0}{\eps}\bigg),\quad
u_t(0,r)=(-c(u(0,r))+\eps)u_r(0,r),
\end{align}
where $\eps>0$ is a constant. Then there exists a small positive constant $\eps_0$ $($the number $\eps_0$ is given in \eqref{2.20}$)$ such that $T<\fr{r_0-\eps}{c_1}$ for any $\eps<\eps_0$.
\end{thm}

The approach of proving Theorem \ref{thm} is inspired by Glassey, Hunter and Zheng \cite{Glassey1} for studying the singularity formation of the one-dimensional equation \eqref{1.4}. We derive an energy equation of smooth solutions for the spherically-symmetric equation \eqref{1.6} which is the key point in this paper. Based on the energy equation, we use the method of characteristics to prove that singularity formation occurs before the effect of the geometric singularity need to be considered. Due to the arbitrary smallness of the initial energy, we know that this is a cusp-type singularity.
Our result is contrasted with the result of Lindblad \cite{Lindblad} for studying the spherically-symmetric nonlinear wave equation
\begin{align}\label{1.10}
u_{tt}-c^2(u)\triangle u=0.
\end{align}
In \cite{Lindblad}, the author showed that there exists a global smooth solution to \eqref{1.10} with smooth, small, and spherically symmetric initial data in $\mathbb{R}^3$. It is clear that, compared with \eqref{1.10}, equation \eqref{1.5} contains a lower-order term $cc'|\nabla u|^2$, which is responsible for the blow-up in the derivatives of $u$. For more related singular examples, one can refer to \cite{Chen-Zheng, CHL, Chen-Huang-Liu, Song} and references therein. We present the detailed proof of Theorem \ref{thm} in the next section.

\section{The proof of Theorem \ref{thm}}\label{S2}

The proof of Theorem \ref{thm} is divided into five steps. In Step 1, we derive the equations in terms of the Riemann variables and the energy equation. In Step 2, we establish the estimate of the energy. Step 3 is focused on establishing the estimate of Riemann variables in the characteristic triangle region. We control the sign of $c'(u)$ in Step 4. Finally, we show that the breakdown of smooth solutions occurs before $t=(r_0-\eps)/c_1$ in Step 5.

\textbf{Step 1.} We introduce the following variables
\begin{align}\label{2.1}
R:=r^\alpha(u_t+c(u)u_r),\quad S:=r^\alpha(u_t-c(u)u_r).
\end{align}
Then, for smooth solutions, there hold by \eqref{1.6}
\begin{align}\label{2.3}
\left\{
\begin{array}{l}
R_t-cR_r=\fr{c'}{4cr^\alpha}(R^2-S^2)-\fr{\alpha c}{r}S, \\
S_t+cS_r=\fr{c'}{4cr^\alpha}(S^2-R^2)+\fr{\alpha c}{r}R.
\end{array}
\right.
\end{align}
We multiply the first equation in \eqref{2.3} by $R$ and the second one by $S$ and then add the resulting equations to obtain the key energy equation
\begin{align}\label{2.4}
(R^2+S^2)_t+(c(u)(S^2-R^2))_r=0.
\end{align}
We point out that the above energy equation may have more
applications, for example, to establish the existence of energy-conservative weak solutions.

\textbf{Step 2.} Corresponding to \eqref{1.9}, we have
\begin{align}\label{2.5}
\begin{array}{l}
u_r(0,r)=\phi'\bigg(\dps\fr{r-r_0}{\eps}\bigg),\quad R(0,r)=\eps r^\alpha u_r(0,r),\\
S(0,r)=(-2c(u(0,r))+\eps)r^\alpha u_r(0,r),
\end{array}
\quad \forall\ r>0.
\end{align}
The small constant $\eps$ can be first chosen $\eps<\min\{c_0,r_0/2\}$. From \eqref{2.5} and \eqref{1.8} we see that
\begin{align}\label{2.6}
R(0,r)=S(0,r)=0,\quad \forall\ r\in[0,r_0-\eps)\cup(r_0+\eps,\infty).
\end{align}
Denote the energy function $E(t)$ as follows
\begin{align}\label{2.7}
E(t):=\int_{0}^\infty [R^2(t,r)+S^2(t,r)]\ {\rm d}r.
\end{align}
Thus from \eqref{2.6} and \eqref{1.7}, we use the energy equation \eqref{2.4} to estimate the energy $E(t)$ for $t\leq (r_0-\eps)/c_1$
\begin{align}\label{2.8}
E(t)=E(0)&=\int_{0}^\infty [R^2(0,r)+S^2(0,r)]\ {\rm d}r \nonumber \\
&= \int_{r_0-\eps}^{r_0+\eps} [\eps^2+(-2c(u(0,r))+\eps)^2]r^{2\alpha}\bigg(\phi'\bigg(\fr{r-r_0}{\eps}\bigg)\bigg)^2\ {\rm d}r
\nonumber \\
&\leq \eps[\eps^2+(2c_1+\eps)^2]r_{0}^{2\alpha}\int_{-1}^1 (\phi'(z))^2\ {\rm d}z\leq Kr_{0}^{2\alpha}\eps
\end{align}
for some positive constant $K$.

\textbf{Step 3.} Let $r_2>r_1>0$ be two constants. We define the positive characteristic curve $r_+(t)$ (or $t_+(r)$) from $(0,r_1)$ and negative characteristic curve $r_-(t)$ (or $t_-(r)$) from $(0,r_2)$ by
\begin{align}\label{2.9}
\left\{
\begin{array}{l}
\fr{{\rm d}r_+(t)}{{\rm d}t}=c(u(t,r_+(t))), \\
r_+(0)=r_1,
\end{array}
\right. \quad
\left\{
\begin{array}{l}
\fr{{\rm d}r_-(t)}{{\rm d}t}=-c(u(t,r_-(t))), \\
r_-(0)=r_2.
\end{array}
\right.
\end{align}
\begin{figure}[htbp]
\begin{center}
\includegraphics[scale=0.55]{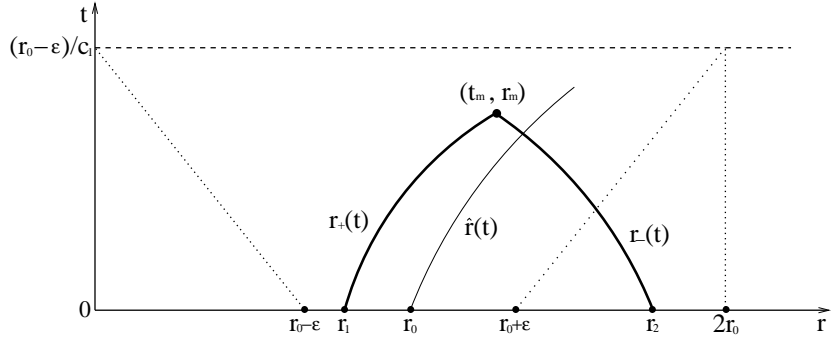}
\caption{\footnotesize The characteristic triangle region.}
\end{center}
\end{figure}
Thanks to \eqref{1.7}, we see that, if
\begin{align}\label{2.10}
r_2-r_1<\fr{2c_0(r_0-\eps)}{c_1},
\end{align}
then the curves $r_+(t)$ and $r_-(t)$ will interact at some point $(t_m,r_m)$ satisfying $t_m< \fr{r_0-\eps}{c_1}$ and $r_1<r_m<r_2$, see Fig. 1. We now consider the characteristic triangle region bounded by $r_\pm(t)$ and the segment $[r_1,r_2]$ and integrate \eqref{2.4} on this characteristic triangle to get by using the divergence theorem
\begin{align}\label{2.11}
\int_{r_1}^{r_m}R^2(t_+(r),r)\ {\rm d}r +\int_{r_m}^{r_2}S^2(t_-(r),r)\ {\rm d}r =\fr{1}{2}\int_{r_1}^{r_2}[R^2(0,r)+S^2(0,r)]\ {\rm d}r,
\end{align}
which combined with \eqref{2.8} gives
\begin{align}\label{2.12}
\int_{r_1}^{r_m}R^2(t_+(r),r)\ {\rm d}r +\int_{r_m}^{r_2}S^2(t_-(r),r)\ {\rm d}r \leq Kr_{0}^{2\alpha}\eps.
\end{align}

\textbf{Step 4.} Denote $r=\hat{r}(t)$ with $t<(r_0-\eps)/c_1$ (or $t=\hat{t}(r)$) the positive characteristic curve from $(0,r_0)$, see Fig. 1. Then the definition of $\hat{r}(t)$ is given by
\begin{align}\label{2.13}
\fr{{\rm d}\hat{r}(t)}{{\rm d}t}=c(u(t,\hat{r}(t))), \quad
\hat{r}(0)=r_0.
\end{align}
We next control the sign of $c'(u)$ on the curve $r=\hat{r}(t)$ by choosing $\eps$ small enough. Recalling the definition of $R$ in \eqref{2.1} yields
\begin{align}\label{2.14}
\fr{{\rm d}u(t,\hat{r}(t))}{{\rm d}t}=\fr{R(t,\hat{r}(t))}{\hat{r}^\alpha(t)}.
\end{align}
We integrate \eqref{2.14} from $0$ to $t<(r_0-\eps)/c_1$ and employ \eqref{2.12} to acquire
\begin{align}\label{2.15}
|u(t,\hat{r}(t))-u(0,r_0)|&=\bigg|\int_{0}^t\fr{R(\tau,\hat{r}(\tau))}{\hat{r}^\alpha(\tau)}\ {\rm d}\tau\bigg| \leq \fr{1}{r_{0}^\alpha}\int_{0}^t|R(\tau,\hat{r}(\tau))|\ {\rm d}\tau \nonumber \\
&\leq \fr{\sqrt{r_0-\eps}}{r_{0}^\alpha\sqrt{c_1}}\cdot\bigg(\int_{0}^tR^2(\tau,\hat{r}(\tau))\ {\rm d}\tau\bigg)^{\fr{1}{2}}  \nonumber \\
&\leq \fr{\sqrt{r_0-\eps}}{r_{0}^\alpha\sqrt{c_1c_0}}\cdot\bigg(\int_{r_0}^rR^2(\hat{t}(r), r)\ {\rm d}r\bigg)^{\fr{1}{2}}\leq \sqrt{\fr{K(r_0-\eps)}{c_0c_1}}\sqrt{\eps}.
\end{align}
According to the $C^2$ regularity assumption of $c(\cdot)$, we can choose $\eps'$ small enough such that for $\eps<\eps'$ there holds on the curve $r=\hat{r}(t)$
\begin{align}\label{2.16}
c'(u(t,\hat{r}(t)))\geq \fr{c'(u(0,r_0))}{2}\geq\fr{c'(u_0)}{4}>0.
\end{align}

\textbf{Step 5.} Now we show that $S(t,\hat{r}(t))$ becomes infinite before time $t=(r_0-\eps)/c_1$. We first recall \eqref{1.8}, \eqref{2.5} and the chosen of $\eps<r_0/2$  to find that
\begin{align}\label{2.17}
S(0,r_0)&=(-2c(u(0,r_0))+\eps)r_{0}^\alpha\phi'(0)> 2 c_0r_{0}^\alpha\max\bigg\{\fr{32c_{1}^22^\alpha}{r_0c_0c'(u_0)},\ \fr{1}{c_0r_{0}^\alpha}\bigg\} \nonumber \\
&>\max\bigg\{\fr{32c_{1}^2(2r_0)^\alpha}{(r_0-\eps)c'(u_0)},\ 2\bigg\},
\end{align}
which means that
\begin{align}\label{2.18}
\fr{1}{S(0,r_0)}<\min\bigg\{\fr{(r_0-\eps)c'(u_0)}{32c_{1}^2(2r_0)^\alpha},\ \fr{1}{2}\bigg\}.
\end{align}
Denote
\begin{align}\label{2.19}
M=\fr{Kc_1r_{0}^\alpha\sqrt{r_0}}{4c_{0}^2}+\fr{\alpha\sqrt{Kc_1}r_{0}^\alpha}{\sqrt{r_0c_0}}.
\end{align}
We choose small constant $\eps_0>0$ satisfying
\begin{align}\label{2.20}
\sqrt{\eps_0}=\min\bigg\{\fr{r_0c'(u_0)}{64Mc_{1}^2(2r_0)^\alpha},\ \fr{1}{2M},\ \sqrt{\eps'}, \ \sqrt{\fr{r_0}{2}},\ \sqrt{c_0}\bigg\}.
\end{align}
Then for $\eps<\eps_0$, we claim that $S(t,\hat{r}(t))$ goes to infinite at some point $t^*<(r_0-\eps)/c_1$. To
this end, we first show that, if $S(t,\hat{r}(t))$ is smooth in $t\in[0,(r_0-\eps)/c_1)$, then $S(t,\hat{r}(t))>1$.
If not, we assume that there exists a number $\tilde{t}\in(0, (r_0-\eps)/c_1)$ such that the function $S(t,\hat{r}(t))\in C^1([0,\tilde{t}])$ satisfying $S(\tilde{t},\hat{r}(\tilde{t}))=1$ and $S(t,\hat{r}(t))>1$ for $t\in[0,\tilde{t})$. By the equation for $S$ in \eqref{2.3}, we derive
\begin{align}\label{2.21}
\fr{\rm d}{{\rm d}t}\bigg(\fr{1}{S}\bigg)\leq & -\fr{c'}{4cr^\alpha}+\fr{1}{S^2}\bigg(\fr{c'}{4cr^\alpha}R^2+\fr{\alpha c}{r}|R|\bigg) \nonumber \\
\leq& -\fr{c'(u_0)}{16c_1(2r_0)^\alpha}+\fr{1}{S^2}\bigg(\fr{c_1}{4c_0r_{0}^\alpha}R^2+\fr{\alpha c_1}{r_0}|R|\bigg).
\end{align}
Integrating \eqref{2.21} from $t=0$ to $t=\tilde{t}$ yields
\begin{align}\label{2.22}
\fr{1}{S(\tilde{t}, \hat{r}(\tilde{t}))}-\fr{1}{S(0,r_0)}\leq \int_{0}^{\tilde{t}}\bigg\{-\fr{c'(u_0)}{16c_1(2r_0)^\alpha} +\fr{1}{S^2}\bigg(\fr{c_1}{4c_0r_{0}^\alpha}R^2+\fr{\alpha c_1}{r_0}|R|\bigg)\bigg\}\ {\rm d}t,
\end{align}
which together with \eqref{2.18} and \eqref{2.12} leads to
\begin{align}\label{2.23}
1=\fr{1}{S(\tilde{t}, \hat{r}(\tilde{t}))}&\leq\fr{1}{S(0,r_0)} + \int_{0}^{\tilde{t}}\fr{1}{S^2}\bigg(\fr{c_1}{4c_0r_{0}^\alpha}R^2+\fr{\alpha c_1}{r_0}|R|\bigg)\ {\rm d}t \nonumber \\
&<\fr{1}{2}+ \fr{c_1}{4c_0r_{0}^\alpha}\int_{r_0}^{\tilde{r}}\fr{R^2}{c}\ {\rm d}r +\fr{\alpha c_1}{r_0}\sqrt{\fr{r_0-\eps}{c_1}}\sqrt{\int_{r_0}^{\tilde{r}}\fr{R^2}{c}\ {\rm d}r} \nonumber \\
&\leq \fr{1}{2} +\bigg(\fr{Kc_1r_{0}^\alpha\sqrt{\eps}}{4c_{0}^2} +\fr{\alpha\sqrt{Kc_1(r_0-\eps)}r_{0}^\alpha}{r_0\sqrt{c_0}}\bigg)\sqrt{\eps}
\nonumber \\
&\leq \fr{1}{2} +\bigg(\fr{Kc_1r_{0}^\alpha\sqrt{r_0}}{4c_{0}^2} +\fr{\alpha\sqrt{Kc_1}r_{0}^\alpha}{\sqrt{r_0 c_0}}\bigg)\sqrt{\eps}\leq \fr{1}{2}+M\sqrt{\eps}<1,
\end{align}
a contradiction. Therefore, we have $S(t,\hat{r}(t))>1$ in $t\in[0,(r_0-\eps)/c_1)$ before the blowup time. We now integrate \eqref{2.21} from $t=0$ to $t<(r_0-\eps)/c_1)$ and apply \eqref{2.18} and \eqref{2.20} to arrive at
\begin{align*}
\fr{1}{S(t, \hat{r}(t))}&\leq \fr{1}{S(0,r_0)}-\fr{c'(u_0)}{16c_1(2r_0)^\alpha}t +\int_{0}^{t}\bigg(\fr{c_1}{4c_0r_{0}^\alpha}R^2+\fr{\alpha c_1}{r_0}|R|\bigg)\ {\rm d}t   \\
&< \fr{(r_0-\eps)c'(u_0)}{32c_{1}^2(2r_0)^\alpha} -\fr{c'(u_0)}{16c_1(2r_0)^\alpha}t +M\sqrt{\eps_0}
 \\
&< \fr{(r_0-\eps)c'(u_0)}{32c_{1}^2(2r_0)^\alpha} -\fr{c'(u_0)}{16c_1(2r_0)^\alpha}t +\fr{r_0c'(u_0)}{64c_{1}^2(2r_0)^\alpha}
  \\
&=\fr{c'(u_0)}{16c_1(2r_0)^\alpha}\bigg(\fr{r_0-\eps}{2c_1}-t+\fr{r_0}{4c_1}\bigg),
\end{align*}
from which we see that there exists a number $t_*$ such that $S(t, \hat{r}(t))\rightarrow\infty$ as $t\rightarrow t_{*}^-$ and the number $t_*$ satisfies
$$
t_*<\fr{r_0-\eps}{2c_1}+\fr{r_0}{4c_1}<\fr{r_0-\eps}{c_1}.
$$
Recalling \eqref{2.1} and noting the boundness of $\hat{r}(t)$, the proof of Theorem \ref{thm} is complete.

\section*{Acknowledgements}

The work of Y. Hu was partially supported by the Zhejiang Provincial Natural Science Foundation (LY17A010019) and National Science Foundation of China (11301128). The work of G. Wang was partially supported by the
Science Key Project of Education Department of Anhui Province (KJ2018A0517).


\end{document}